\documentclass[12pt]{amsart}
\usepackage[dvips]{graphics}
\usepackage[latin1]{inputenc}
\usepackage{amssymb, amsmath, amscd, mathrsfs,marvosym, psfrag}
\input xy
\xyoption{all}
\DeclareMathAlphabet{\mathpzc}{OT1}{pzc}{m}{it}

\newtheorem{theorem}{Theorem}[section]

\newtheorem{proposition}[theorem]{Proposition}

\newtheorem{conjecture}[theorem]{Conjecture}
\newtheorem{lemma}[theorem]{Lemma}
\newtheorem*{ThM}{Main Theorem}

\theoremstyle{definition}
\newtheorem{definition}[theorem]{Definition}

\theoremstyle{remark}

\newcommand{\CB}{{\mathcal B}}

\newcommand{\CE}{{\mathcal E}}

\newcommand{\CG}{{\mathcal G}}

\newcommand{\CM}{{\mathcal M}}
\newcommand{\CN}{{\mathcal N}}

\newcommand{\CT}{{\mathcal T}}

\newcommand{\CV}{{\mathcal V}}
\newcommand{\CW}{{\mathcal W}}

\newcommand{\CZ}{{\mathcal Z}}

\newcommand{\SB}{{\mathscr B}}

\newcommand{\SL}{{\mathscr L}}
\newcommand{\SM}{{\mathscr M}}
\newcommand{\SN}{{\mathscr N}}

\newcommand{\SZ}{{\mathscr Z}}

\newcommand{\fg}{{{\mathfrak g}}}

\newcommand{\hCW}{{\widehat\CW}}

\newcommand{\hCG}{{\widehat\CG}}

\newcommand{\hCT}{{\widehat\CT}}

\newcommand{\hV}{{\widehat V}}
\newcommand{\hR}{{\widehat R}}
\newcommand{\hX}{{\widehat X}}

\newcommand{\hH}{{\widehat H}}

\newcommand{\DZ}{{\mathbb Z}}

\newcommand{\DN}{{\mathbb N}}
\newcommand{\DQ}{{\mathbb Q}}

\newcommand{\ch}{{\operatorname{char}\, }}

\newcommand{\Hom}{{\operatorname{Hom}}}

\newcommand{\res}{{\operatorname{res}}}

\newcommand{\supp}{{\operatorname{supp}}}
\newcommand{\catmod}{{\operatorname{-mod}}}

\newcommand{\inj}{{\hookrightarrow}}

\newcommand{\Loc}{{\SL}}

\newcommand{\dual}{{\mathsf D}}
\newcommand{\Zd}{{{\mathsf D}{\CZ}}}
\newcommand{\Bd}{{{\mathsf D}{\CB}}}

\newcommand{\grk}{{{\operatorname{\underline{rk}}}}}

\newcommand{\ol}{\overline}

\newcommand{\Sh}{\operatorname{Sh}}

\newcommand{\linie}{{\,\text{---\!\!\!---}\,}}
\newcommand{\llinie}{{\text{---\!\!\!---\!\!\!---}}}

\begin{document}

\pagenumbering{arabic}
\title[]{The multiplicity one case of Lusztig's conjecture}
\author[]{Peter Fiebig}
\address{Mathematisches Institut, Universit{ä}t Freiburg, 79104 Freiburg, Germany}
\email{peter.fiebig@math.uni-freiburg.de}

\begin{abstract} We prove the multiplicity one case of Lusztig's conjecture on the irreducible characters of reductive algebraic groups for all fields with characteristic above the Coxeter number. 
\end{abstract}
\maketitle


\section{Introduction}
Some deep results in representation theory rely on the possibility to interpret representation theoretic structures geometrically and to make use  of the powerful methods of algebraic geometry. One of the most important results applied in such situations is the decomposition theorem of Beilinson, Bernstein, Deligne and Gabber. It assumes, however, the coefficients of the theory to be of characteristic zero. While it is often possible to deduce results for ``almost all'' characteristics from their characteristic zero analogues, the decomposition theorem does not convey the exceptional primes. 

\subsection{Lusztig's conjecture}

In 1980 Lusztig conjectured a formula for the irreducible rational
characters of reductive algebraic groups over a field of positive
characteristic (cf.~ \cite{L1}). In the early 90s,  a combined effort
of several authors that followed a program outlined by Lusztig
culminated in a proof of the conjecture for almost all
characteristics. This means, more precisely, that for a fixed root system the formula is correct if the characteristic of the ground field is larger than a certain bound. This bound, however, is unknown in all but low rank cases. 

The central idea of Lusztig's program was to first prove the 
characteristic zero case of the conjecture using the above mentioned
geometric approach. This involved the decomposition theorem for
sheaves on  affine Schubert varieties with complex coefficients, Kac--Moody algebras and quantum groups. A result of Andersen, Jantzen and Soergel then provided the means to deduce Lusztig's original conjecture from its quantum version for almost all characteristics. 

\subsection{Lusztig's conjecture as a moment graph problem}
The papers \cite{ModRep} and \cite{CharIrr} provide a new approach
towards Lusztig's conjecture which is largely independent of the
program outlined above. The main result is a direct link between
certain sheaves on affine Schubert varieties with coefficients in a
field $k$ of positive characteristic and  representations of the Lie
algebra over $k$ associated to the dual root system. An intermediate
step in the construction of this link involves the Braden-MacPherson
sheaves on an associated moment graph. 

It is well-known (and proven in \cite{CharIrr}) that Lusztig's conjecture is equivalent to a certain
conjecture on the Jordan--H\"older multiplicities of baby Verma
modules for the Lie algebra (cf.~Conjecture \ref{conj-multipl}). In \cite{ModRep} we
showed that the latter conjecture is implied by a  multiplicity conjecture on the ranks of the stalks of the
Braden-MacPherson sheaves (cf.~ Conjecture \ref{conj-MC}). Apart from
giving a new proof for almost all characteristics, this result is used
in  \cite{CharBo} to give a specific upper bound  for the set of exceptional primes. This bound, however, is huge. 

\subsection{The multiplicity one case}
The purpose of the present article is to prove the multiplicity one
case of the multiplicity conjecture mentioned above. Before we
formulate our result we introduce some notation. Let $R$ be a reduced,
irreducible and finite root system and let $X$ be its weight
lattice. Let $\fg$ be the Lie algebra associated to $R$ over the field
$k$. After fixing a set of positive roots we have for each $\lambda\in
X$  the baby Verma module $\Delta(\lambda)$ and its simple quotient
$L(\lambda)$. We denote by  $\hCW$ the affine Weyl group. It acts by affine transformations on the lattice $X$. Let $\rho\in X$ be the half-sum of positive roots. Then Lusztig's conjecture is equivalent to the following multiplicity conjecture:

\begin{conjecture}\label{conj-multipl} Suppose that the characteristic of $k$ is bigger than the Coxeter number of $R$. For all $x,w\in\hCW$ we have 
$$
[\Delta(x(\rho)-\rho):L(w(\rho)-\rho)]=p_{x,w}(1).
$$
\end{conjecture}

Here  $p_{x,w}$ denotes the periodic polynomial defined by Lusztig (\cite{L2}, cf.~  \cite{CharIrr} for notation and normalization). The main result of this article is the following.

\begin{theorem} Suppose that the characteristic of $k$ is bigger than the Coxeter number of $R$. For all $x,w\in\hCW$ we  have
$$
[\Delta(x(\rho)-\rho):L(w(\rho)-\rho)]=1\text{ if and only if } p_{x,w}(1)=1.
$$
\end{theorem}

\subsection{The smooth locus of a moment graph}
The above theorem is an application of a more general result on the
smooth locus of a moment graph. Let us quickly recall the main
notions and results on sheaves on moment graphs. 

A moment graph $\CG$ is a graph, whose edges are labelled by non-trivial elements of a vector space $V$, and that is endowed with a partial order on its set of vertices such that two vertices are comparable if they lie on a common edge.

There are two naturally defined sheaves on $\CG$ (cf.~ Section
\ref{sheaves}): the first  is the structure sheaf $\SZ$, and the
second is the  Braden-MacPherson sheaf $\SB$. Taking global sections
gives the cohomology $\CZ$ of the graph, which is a graded algebra over the symmetric algebra $S=S(V)$, and its intersection cohomology $\CB$, which is  a graded $\CZ$-module. For each vertex $x$ of $\CG$ there are stalks $\CZ^x$ and $\CB^x$ and, by definition, $\CZ^x=S$ and $\CB^x$ is isomorphic to a graded free $S$-module of finite rank. More precisely, $\CB^x\cong S\oplus\bigoplus S\{k_i\}$ for some $k_i<0$. We define the {\em  smooth locus} of $\CG$ to be the set of vertices $x$ with $\CB^x\cong S$. 

There are far too many moment graphs to hope for a general and satisfactory theory. The moment graphs that are related to Lusztig's conjecture have an additional,  distinguishing feature: their intersection cohomology is self-dual up to a degree shift.  Here is our main result: 

\begin{ThM} Assume that $\CB$ is self-dual of degree $l$, i.e.~ $\Hom^\bullet_S(\CB,S)\cong \CB\{l\}$. Let $x\in\CV$. Then $\CB^x\cong S$ if and only if  for all $y\geq x$ there are exactly $l$ edges of $\CG$ containing $y$. 
\end{ThM}

Note that, once it is known that $\CB$ is self-dual, the theorem only refers to the underlying abstract graph structure and the partial order, so ignores the labelling on $\CG$. In the case of moment graphs of Coxeter systems, this phenomenon is related to the {\em combinatorial invariance conjecture} for Kazhdan--Lusztig polynomials (cf.~ \cite{BB}).

\subsection{Contents}

The paper is organized as follows. In Section \ref{SHM} we introduce sheaves on   moment graphs and the  structure algebra $\CZ$. In Section \ref{Zmod} we discuss modules over $\CZ$ and their relation to sheaves. In Section \ref{IC} we define the intersection cohomology $\CB$ and recall some of its properties. In Section \ref{smloc} we study the smooth locus of a moment graph and prove  the main theorem. Section \ref{Lconj} contains the application to Lusztig's conjecture.

\section{Sheaves on moment graphs}\label{SHM}

In this section we review the theory of sheaves on moment graphs. References for the following are \cite{BMP} and \cite{Fie1}.
\subsection{Moment graphs}
Let $k$ be a field and $V$ a finite dimensional $k$-vector space. 
A {\em moment graph} $\CG=(\CV,\CE,\alpha,\leq)$ over $V$ is the following data:
\begin{itemize}
\item a graph $(\CV,\CE)$ given by a set $\CV$ of vertices and a set $\CE$ of edges,
\item a labelling $\alpha\colon \CE\to V\setminus\{0\}$ of each edge  by a non-zero element of $V$,
\item a partial order ``$\leq$'' on $\CV$ such that two vertices are comparable if they are connected by an edge.
\end{itemize}

The partial order induces a direction on each edge. We write $E\colon x\to y$ to denote an edge $E$ with endpoints $x,y$ such that $x<y$, or $E\colon x\stackrel{\alpha}\longrightarrow y$ if, in addition,  we want to indicate the label $\alpha=\alpha(E)$. If we want to ignore the direction we write $E\colon x\linie y$ or  $E\colon x\stackrel{\alpha}\llinie y$.

We assume that there are no double edges in $(\CV,\CE)$, i.e.~ that two vertices are connected by at most one edge, and that the partial order is the same as the partial order generated by the directed edges (i.e.~ two comparable vertices are linked by a directed path in $\CG$). 
Most of the constructions that we describe in this paper suppose that $\CG$ is finite and that there is a highest element in the partial order, i.e.~ an element $w\in\CV$ such that $x\leq w$ for all vertices $x\in\CV$. Moreover, we assume that the moment graph has the following property:

\begin{definition} A moment graph $\CG$ is called a {\em GKM-graph} if the labels of any two distinct edges containing a common vertex are linearly independent. 
\end{definition}

For our approach the GKM-property is essential.
In Section \ref{Lconj} we associate a moment graph to an affine root system over an arbitrary field $k$. The GKM-property then defines the only restriction on the characteristic of $k$ that we need. An important feature that distinguishes the moment graphs that we are ultimately interested in is the {\em self-duality of the intersection cohomology} (cf.~ Section \ref{IC}). 

\subsection{The structure algebra of a moment graph}

Let $\CG$ be a finite moment graph over $V$. Let $S=S_k(V)$ be the symmetric algebra over $V$, $\DZ$-graded such that $\deg V=1$. (Note that it would be much more customary to define $\deg V=2$ (cf.~ \cite{Soe04,Fie2}), but the normation above eases the exposition in this article considerably.) All $S$-modules that we consider in the following are $\DZ$-graded, and all homomorphisms $f\colon M\to N$ between graded $S$-modules are homogeneous of degree zero. Let us define the {\em structure algebra} of $\CG$:
$$
\CZ=\CZ(\CG):=\left\{\left. (z_x)\in\bigoplus_{x\in\CV}S\, \right|\,\begin{matrix} z_x\equiv z_{y}\mod \alpha \\ \text{for all edges $x\stackrel{\alpha}\llinie y$ of $\CG$} \end{matrix}\right\}.
$$
There is a diagonal inclusion $S\inj\CZ$, and coordinatewise
multiplication gives  $\CZ$ the structure of an $S$-algebra. The grading on $S$
induces a grading on $\CZ$. 

Let $\Lambda$ be a full subgraph of $(\CV,\CE)$. In the following we do not distinguish between full subgraphs and subsets of $\CV$. Define 
$$
\CZ_\Lambda:=\left\{(z_x)\in\CZ\mid z_x=0 \text{ for $x\not\in\Lambda$}\right\},
$$ 
the subalgebra of elements supported on $\Lambda$, and 
$$
\CZ^\Lambda:=\CZ/\CZ_{\CV\setminus\Lambda}
.
$$
Note that $\CZ^\Lambda$ is the image of the composition $\CZ\subset\bigoplus_{x\in\CV} S\stackrel{p}\to \bigoplus_{x\in\Lambda} S$,
where $p$ is the projection along the decomposition.
There is a natural inclusion $\CZ_\Lambda\subset\CZ^\Lambda$, and for $\Lambda^\prime\subset\Lambda$ there is a natural inclusion $\CZ_{\Lambda^\prime}\subset\CZ_\Lambda$ and a natural surjection $\CZ^{\Lambda}\to\CZ^{\Lambda^\prime}$. We write $\CZ^{x,y,\dots,z}$ and $\CZ_{x,y,\dots,z}$ instead of $\CZ^{\{x,y\dots,z\}}$ and $\CZ_{\{x,y,\dots,z\}}$, resp. 

Let $Q$ be the quotient field of $S$ and set $\CZ_Q:=\CZ\otimes_S Q$. The inclusion $\CZ\inj\bigoplus_{x\in\CV} S$ induces an inclusion $\CZ_Q\inj\bigoplus_{x\in\CV}Q$. Recall that we can consider $\CZ_x$ as a subset of $\CZ^x$.
\begin{lemma}\label{strucZ} 
\begin{enumerate} 
\item For each vertex $x\in\CV$ we have $\CZ^x=S$ and $\CZ_x=\alpha_1\cdots\alpha_n\cdot\CZ^x$, where $\alpha_1,\dots,\alpha_n$ are the labels of the edges containing $x$.
\item The natural inclusion $\CZ_Q\inj\bigoplus_{x\in\CV} Q$ is a bijection.
\end{enumerate}
\end{lemma}

\begin{proof} The subspace $\CZ^x\subset S$ is an $S$-submodule that contains the generator $1$ as the image of the constant section $(1_x)\in\CZ$. Hence $\CZ^x=S$. Obviously $\alpha_1\cdots\alpha_k\in\CZ_x\subset S$. Conversely, each $f\in\CZ_x$ must be divisible in $S$ by each of the $\alpha_i$. Recall that we assume that $\CG$ is a GKM-graph. So the $\alpha_i$ are pairwise linearly independent, hence $f$ must be divisible by $\alpha_1\cdots\alpha_k$, so $\CZ_x=\alpha_1\cdots\alpha_k\cdot\CZ^x$. Hence we proved (1).  Claim (2) follows easily from the second statement in (1).
\end{proof}

\subsection{Sheaves on a moment graph}\label{sheaves}
A {\em sheaf on $\CG$} is given by the data $\SM=(\SM^x,\SM^E,\rho_{x,E})$, where
\begin{itemize} 
\item $\SM^x$ is a finitely generated, graded and torsion free $S$-module for any vertex $x\in\CV$,
\item $\SM^E$ is a finitely generated, graded $S$-module with $\alpha\cdot M=(0)$ for any edge $E$ with label $\alpha$, and
\item for each vertex $x$ of $E$, $\rho_{x,E}\colon\SM^x\to\SM^E$ is a homomorphism of graded $S$-modules. 
\end{itemize}
A {\em morphism of sheaves}, $f\colon\SM\to\SN$, consists of homomorphisms $f^x\colon \SM^x\to\SN^x$ and $f^E\colon \SM^E\to\SN^E$ of $S$-modules for each vertex $x$ and each edge $E$, such that $\rho_{x,E}^{\SN}\circ f^x=f^{E}\circ \rho^{\SM}_{x,E}$, whenever $x$ is contained in $E$.

We define the {\em structure sheaf} $\SZ$ on $\CG$ by $\SZ^x=S$ for any vertex $x$, $\SZ^E=S/\alpha\cdot S$ for any edge $E$ labelled by $\alpha$, and $\rho_{x,E}\colon S\to S/ \alpha\cdot S$ the canonical map for any vertex $x$ of $E$.

For a sheaf $\SM$ on $\CG$  let
$$
\Gamma(\SM)=\left\{\left.(m_x)\in\bigoplus_{x\in\CV} \SM^x\right|\,\begin{matrix} \rho_{x,E}(m_x)=\rho_{y,E}(m_y) \\ \text{ for all edges $E\colon x\linie y$ }\end{matrix}\right \}
$$
be the space of {\em global sections} of $\SM$. By definition, $\CZ=\Gamma(\SZ)$. More generally, each $\Gamma(\SM)$ is a $\CZ$-module under coordinatewise multiplication. For any subset $\Lambda\subset\CV$ we define
$$
\Gamma(\Lambda,\SM)=\left\{\left.(m_x)\in\bigoplus_{x\in\Lambda} \SM^x\right|\,\begin{matrix} \rho_{x,E}(m_x)=\rho_{y,E}(m_y)\text{ for all} \\ \text{edges $E\colon x\linie y$ with $x,y\in\Lambda$}\end{matrix}\right \}.
$$

\subsection{Local and global viewpoints}
 While  a sheaf  $\SM$ is a local object, the corresponding $\CZ$-module $\Gamma(\SM)$ is a global object. In \cite{Fie1} we constructed a left inverse of the functor $\Gamma$, the {\em localization functor} $\Loc$. We defined the categories $\Sh^{glob}$ of sheaves that are {\em defined by global sections}, and $\CZ\catmod^{loc}$ of $\CZ$-modules that are {\em determined by local relations}, as the images of the respective functors. Moreover, we showed that $\Gamma$ and $\Loc$ are mutually inverse equivalences between $\Sh^{glob}$ and $\CZ\catmod^{loc}$. 

Hence we can choose a local or a global viewpoint on such objects. This is reflected, for example, in the two essentially different constructions of the intersection cohomology for moment graphs associated to Coxeter systems (cf.~ Section \ref{Lconj}). In order to apply our main result in this situation, it is important to have both constructions at hand.

\section{$\CZ$-modules}\label{Zmod}

Let $M$ be a graded $\CZ$-module. We denote by $M_{\{n\}}\subset M$ its homogeneous component of degree $n$, so $M=\bigoplus_{n\in\DZ} M_{\{n\}}$. For $m\in\DZ$ let $M\{m\}$ be the graded $\CZ$-module obtained from $M$ by shifting the grading, i.e.~ such that $M\{m\}_{\{n\}}=M_{\{m+n\}}$. Let $\CZ\catmod^f$ be the category of graded $\CZ$-modules that are finitely generated over $S$ and torsion free over $S$.

\subsection{Decomposition of $\CZ$-modules at the generic point}
Recall that $Q$ is the quotient field of $S$. For each $S$-module $M$ set $M_Q:=M\otimes_S Q$. If $M\in\CZ\catmod^f$, then $M_Q$ is a $\CZ_Q$-module that is finitely generated as a $Q$-vectorspace. By Lemma \ref{strucZ} we have $\CZ_Q=\bigoplus_{x\in\CV}Q$. The idempotents in $\bigoplus_{x\in\CV}Q$ hence give a canonical decomposition 
$$
M_Q=\bigoplus_{x\in\CV} M_{Q,x}.
$$
Since $M$ is supposed to be torsion free over $S$, the natural map
$M\to M_Q$, $m\mapsto m\otimes 1$, is an inclusion, so we can view each $m\in M$ as a tuple $m=(m_x)$ with $m_x\in M_{Q,x}$. Let $\supp(m):=\{x\in\CV\mid m_x\neq 0\}$ and $\supp(M):=\bigcup_{m\in M} \supp(m)$ be the {\em  support} of $m$ and $M$, resp.  

For a subset $\Lambda\subset\CV$ we define
$$
M_\Lambda:= M\cap \bigoplus_{x\in\Lambda} M_{Q,x}, \quad M^\Lambda:=M/M_{\CV\setminus \Lambda}.
$$ 
$M_\Lambda$ and $M^\Lambda$ are again $\CZ$-modules (even $\CZ^\Lambda$-modules) and as such are objects in $\CZ\catmod^f$. We write  $M_{x,y,\dots,z}$ and $M^{x,y,\dots,z}$ instead of $M_{\{x,y\dots,z\}}$ and $M^{\{x,y,\dots,z\}}$, resp. Each homomorphism $f\colon M\to N$ induces homomorphisms $f_\Lambda\colon M_\Lambda\to N_\Lambda$ and $f^\Lambda\colon M^\Lambda\to N^\Lambda$ for each $\Lambda\subset\CV$.

For a vertex $x\in\CV$ define $\CM(x)$ as the graded free $S$-module
of rank 1, generated in degree zero, on which $(z_y)\in \CZ$ acts by
multiplication with $z_x\in S$. It is the space of global sections of
the sheaf $\SM(x)$ that is defined by $\SM(x)^x=S$ and $\SM(x)^y=0$ if
$y\neq x$, and $\SM(x)^E=0$ for all edges $E$. Note that all
homomorphisms $\CM(x)\to\CN$ and $\CN\to\CM(x)$ in $\CZ\catmod^f$
factor over $\CN_x\to\CN$ and $\CN\to\CN^x$, resp. 

\subsection{The dual of the structure algebra}\label{dualofZ}

For any $M\in\CZ\catmod^f$ set 
$$
\dual M:=\Hom^\bullet_S(M,S)=\bigoplus_{n\in\DZ}\Hom^{n}_S(M,S),
$$
where $\Hom^ {n}_S(M,S):=\Hom_S(M,S\{n\})$. Then $\dual M$ is a graded $\CZ$-module with the action given by $(z.f)(m)=f(z.m)$ for $z\in\CZ$, $f\in\dual M$ and $m\in M$. 

\begin{definition} Let $l\in\DZ$. We say that $M\in\CZ\catmod^f$ is {\em self-dual of degree $l$} if there is an isomorphism $\dual M\cong M\{l\}$ of graded $\CZ$-modules. 
\end{definition}

In this section we are particularly  interested in the dual of the
structure algebra itself. Denote by $\Hom_S(\CZ,S)$ the space of
ungraded homomorphisms. Note that this is the same as the space
$\Hom_S^\bullet(\CZ,S)$ with the grading forgotten, as $\CZ$ is a
finitely generated $S$-module. Consider the following:
$$
\Zd=\Hom_S(\CZ,S)\subset\Hom_Q(\CZ_Q,Q)=\bigoplus_{x\in\CV} Q.
$$
Here we identify $\Hom_S(\CZ,S)$ with the subset in $\Hom_Q(\CZ_Q,Q)$ consisting of all $\phi$ such that $\phi(\CZ)\subset S$. Lemma \ref{strucZ} gives the identity $\Hom_Q(\CZ_Q,Q)=\bigoplus_{x\in\CV} Q$. We get a canonical identification
$$
\Zd=\left\{(\phi_x)\in\bigoplus_{x\in\CV} Q\mid \sum_{x\in\CV} \phi_x\cdot z_x\in S\text{ for all $(z_x)\in\CZ$}\right\}.\eqno{(\ast)}
$$
\begin{lemma}\label{strucDZ} Choose $x\in\CV$ and denote by $\alpha_1,\dots,\alpha_n$ the labels of the edges containing $x$. Then we have 
$$
S=(\Zd)_x\subset(\Zd)^x\subset\alpha_1^{-1}\cdots\alpha_n^{-1}\cdot S\subset Q.
$$ 
In particular, $(\Zd)^x$ lives in degrees $\geq -n$. 
\end{lemma}
 
\begin{proof} Let $\phi=(\phi_y)\in(\Zd)_x$, so $\phi_y=0$ for $y\neq x$. Then evaluation at the constant section $1_\CZ$ shows $\phi_x\in S$. Conversely, each $\phi=(\phi_x)$ with $\phi_y=0$ for $y\neq x$ and $\phi_x\in S$ defines an element in $(\Zd)_x$ by $(\ast)$. This shows $S=(\Zd)_x$. 

Let $\phi=(\phi_y)\in\Zd$. By Lemma \ref{strucZ}, $\CZ_x=\alpha_1\cdots\alpha_n\cdot S$, so evaluating $\phi$ on $\CZ_x$ shows $\phi_x\in\alpha_1^{-1}\cdots\alpha_n^{-1}\cdot S$, hence $(\Zd)^x\subset\alpha_1^{-1}\cdots\alpha_n^{-1}\cdot S$.
\end{proof}

\subsection{The completion of a $\CZ$-module} In  \cite{Fie1} we introduced a localization functor $\Loc$ that associates a sheaf on $\CG$ to each $\CZ$-module and that is left adjoint to the functor $\Gamma$. We define the {\em completion} of a $\CZ$-module $M$ as the $\CZ$-module $\widehat{M}=(\Gamma\circ\Loc)(M)$. The adjunction morphism gives a natural map $M\to\widehat M$. In this section we give an explicit description of $M\to \widehat M$.

Let $E\colon x\stackrel{\alpha}\llinie y$ be an edge of $\CG$ and define the {\em local structure algebra} at $E$:
$$
\CZ(E):=\{(z_x,z_y)\in S\oplus S\mid z_x\equiv z_y \mod\alpha\}.
$$
Let 
$M\in\CZ\catmod^f$ and consider the inclusion $M^{x,y}\subset M^x\oplus M^y$. The local structure algebra acts on $M^x\oplus M^y$ by coordinatewise multiplication. Let 
$$
M(E):=\CZ(E)\cdot M^{x,y}\subset M^x\oplus M^y
$$ 
be the submodule generated by $M^{x,y}$. (Note that the map $\CZ\to\CZ(E)$, $z\mapsto (z_x,z_y)$ need not be surjective, so $M^{x,y}$ need not be a $\CZ(E)$-module.) For $M,N\in\CZ\catmod^f$ and a homomorphism $f\colon M\to N$ we have an induced homomorphism $f(E)\colon M(E)\to N(E)$.  

Define 
$$
\widehat{M}:=\left\{\left. (m_x)\in\bigoplus_{x\in\CV}M^x\, \right|\,\begin{matrix} (m_x,m_y)\in M(E) \text{ for } \\ \text{all edges $E\colon x\linie y$} \end{matrix}\right\}.
$$
$\widehat{M}$ is again a $\CZ$-module that we call the {\em completion} of $M$. Note that there is an obvious inclusion $M\inj \widehat{M}$ of $\CZ$-modules. We say that $M$ is {\em complete} if this inclusion is a bijection.

The next lemma follows easily from the definitions.
\begin{lemma} For each vertex $x\in\CV$ we have $(\widehat{M})^x=M^x$. For each edge $E\colon x\linie y$ we have $\widehat{M}(E)=M(E)$. Hence $\widehat{\widehat{M}}=\widehat{M}$.
\end{lemma}

The next statement is needed in Section \ref{propB}.
\begin{lemma}\label{loccom} Let $M=\Gamma(\SM)$ for a sheaf $\SM$ on $\CG$. Then $M$ is complete.
\end{lemma}
\begin{proof} Let $E\colon x\linie y$ be an edge and $(m_x,m_y)\in \Gamma(\SM)(E)$. Then $\rho_{x,E}(m_x)=\rho_{y,E}(m_y)$, since $\Gamma(\SM)(E)=\CZ(E)\cdot\Gamma(\SM)^{x,y}$. Hence $(m_x)\in\widehat{\Gamma(\SM)}$ implies $(m_x)\in\Gamma(\SM)$. 
\end{proof}
\subsection{Flabby objects}
So far we did not make use of the partial order on $\CV$ which is one of the main ingredients of a moment graph. We now use it to define a topology on $\CV$:
\begin{definition}
Let $\Omega\subset\CV$ be a subset. It is called {\em open} if it is upwardly closed, i.e.~ if $\Omega=\bigcup_{x\in\Omega}\{\geq x\}$, where $\{\geq x\}:=\{y\in\CV\mid y\geq x\}$. 
\end{definition}

The next definition is central to our approach.

\begin{definition} Let $M\in\CZ\catmod^f$. $M$ is called {\em flabby}, if for any open $\Omega\subset\CV$ the module $M^\Omega$ is complete. 
\end{definition}

Let $\SM=\Loc(M)$ be the sheaf associated to $M$. In \cite[Proposition 4.2]{Fie1} it is shown that $M$ is flabby if and only if $M=\Gamma(\SM)$ and the natural map $\Gamma(\SM)\to\Gamma(\Omega,\SM)$ is surjective for all open $\Omega\subset\CV$. This justifies the term ``flabby''.

Note that structure algebra $\CZ$ itself may fail to be flabby.  

\section{The intersection cohomology of $\CG$}\label{IC}
In this section we construct, following \cite{BMP}, a $\CZ$-module $\CB$ which can be thought of as correcting $\CZ$'s failure of being flabby. It will be given as the space of global sections of a sheaf $\SB$ on $\CG$. We assume that $\CG$ is finite and that there is a highest vertex $w$ in the partial order, i.e.~ such that $x\leq w$ for all $x\in\CV$.  

\subsection{The Braden--MacPherson construction} We want to construct a flab­by sheaf on $\CG$, i.e.~ a sheaf with the property that all sections that are defined on an open subset extend to global sections.  
Let $\SM$ be a sheaf on $\CG$ and choose $x\in\CV$. Set $\{>x\}=\{y\in\CV\mid y>x\}$ and define $\{\geq x\}$ similarly. Consider the restriction map
$$
\Gamma(\{\geq x\},\SM)\to\Gamma(\{>x\},\SM).\eqno{(\ast\ast)}
$$

Suppose $(m_y)\in\Gamma(\{>x\},\SM)$ is a section. In order to extend it to the vertex $x$ we have to find $m_x\in\SM^x$ such that $\rho_{x,E}(m_x)=\rho_{y,E}(m_y)$ for all $y>x$ that are connected to $x$ by an edge. This leads us to the following definitions.

We define $\CV_{\delta x}\subset\CV$ as the set of vertices $y>x$ that are connected to $x$ by an edge. Accordingly, let $\CE_{\delta x}=\{E\mid E\colon x\to y, y\in\CV_{\delta x}\}\subset\CE$ be the set of the corresponding edges.
We define $\SM^{\delta x}\subset\bigoplus_{E\in\CE_{\delta x}}\SM^E$ as the image of  the composition
$$
\Gamma(\{>x\},\SM)\subset\bigoplus_{y>x}\SM^y \stackrel{p}\to\bigoplus_{y\in\CV_{\delta x}} \SM^y\stackrel{\rho}\to\bigoplus_{E\in\CE_{\delta x}}\SM^E,
$$
where $p$ is the projection along the decomposition and $\rho=\bigoplus_{y\in\CV_{\delta x}}\rho_{y,E}$. 
Finally, we define  
$$
\rho_{x,\delta x}:=(\rho_{x,E})_{E\in\CE_{\delta x}}^T\colon \SM^x\to\bigoplus_{E\in\CE_{\delta x}}\SM^E.
$$
So the restriction map $(\ast\ast)$  is surjective if and only if $\SM^{\delta x}$ is contained in the image of $\rho_{x,\delta x}$.

\begin{theorem}[\cite{BMP}]\label{conB} There is an up to isomorphism unique sheaf $\SB$ on $\CG$ with the following properties:
\begin{enumerate}
\item $\SB^w\cong S$.
\item If $x\stackrel{\alpha}\to y$ is a directed edge, then the map $\rho_{y,E}\colon \SB^y\to\SB^E$ is surjective with kernel $\alpha\cdot\SB^y$.
\item For any $x\in\CV$, the image of $\rho_{x,\delta x}$ is $\SB^{\delta x}$, and  $\rho_{x,\delta x}\colon \SB^x\to\SB^{\delta x}$ is a projective cover in the category of graded $S$-modules. 
\end{enumerate}
\end{theorem}
\begin{proof} We can actually use the stated properties to define $\SB$ inductively, as follows. We start with setting $\SB^w=S$. Once we defined $\SB^y$, we set $\SB^E=\SB^y/ \alpha\cdot\SB^y$ for any edge $E\colon x\stackrel{\alpha}\to y$. For $x\in\CV$ let  $\CG_{>x}$ be the full subgraph given by the set of vertices $\{>x\}$. Once we defined $\SB$ on $\CG_{>x}$, the $S$-module $\SB^{\delta x}\subset\bigoplus_{E\in\CE_{\delta x}}\SB^E$ can already be constructed, and we define $\SB^x\to\SB^{\delta x}$ as a projective cover. This also induces the maps $\rho_{x,E}\colon \SB^x\to\SB^E$ as the components of the projective cover map.

We now prove the uniqueness. If $\SB^\prime$ is another sheaf satisfying the stated properties, we can choose an isomorphism $\SB^w\cong(\SB^\prime)^w$. Once we have an isomorphism $\SB^y\cong(\SB^{\prime})^y$, there is a unique induced isomorphism $\SB^E\cong(\SB^\prime)^E$ for each edge $E\colon x\stackrel{\alpha}\to y$. Once we defined the isomorphism on the full subgraph $\CG_{>x}$, we get an induced isomorphism $\SB^{\delta x}\cong(\SB^\prime)^{\delta x}$ and, finally, a compatible isomorphism between the projective covers $\SB^x$ and $(\SB^\prime)^x$.
\end{proof}

 We call $\SB$ the {\em Braden-MacPherson sheaf}, and $\CB:=\Gamma(\SB)$ the {\em intersection cohomology} of $\CG$. This is justified by the following: In \cite{BMP} a moment graph was associated to certain complex varieties endowed with an action of a complex torus $T$. It was shown that $\CB$ encodes the structure of the $T$-equivariant intersection cohomology of the variety.

From the construction it follows that each local section of $\SB$ extends to a global section. In particular, we have $\CB^x=\SB^x$ for all vertices $x$. We define $\CB^E:=\SB^E$ for an edge $E$, and $\CB^{\delta x}:=\SB^{\delta x}$ for a vertex $x$. 

\subsection{Properties of $\CB$}\label{propB} The following proposition shows that $\CB$ and $\CZ$ coincide if and only if $\CZ$ is flabby.  

\begin{proposition}\label{flabbyB}
\begin{enumerate}
\item For each open subset $\Omega$ of $\CV$ we have $\CB^\Omega=\Gamma(\Omega,\SB)$.
\item  $\CB$ is flabby. 
\item There is an isomorphism $\CZ\cong \CB$ if $\CZ$ is flabby.
\end{enumerate}
\end{proposition}
\begin{proof} Since each local section in $\Gamma(\Omega,\SB)$ (for $\Omega\subset\CV$ open) extends to a global section, we have $\CB^\Omega=\Gamma(\Omega,\SB)$. Hence $\CB^\Omega$ is complete by Lemma \ref{loccom}, thus $\CB$ is flabby. 

In any case we have $\rho_{x,\delta x}(\SZ^x)\subset\SZ^{\delta x}$ for all $x$, since the image of $1$ lies in $\SZ^{\delta x}$. If $\CZ$ is flabby, then $\rho_{x,\delta x}\colon \SZ^x\to\SZ^{\delta x}$ is surjective. Since $\SZ^x=S$, $\rho_{x,\delta x}$  must be a projective cover. It is clear that $\SZ$ also satisfies the two remaining properties listed in Theorem \ref{conB}, so $\SZ\cong \SB$ by uniqueness.
\end{proof}

\begin{lemma}\label{indB} Let $\Omega\subset\CV$ be an open subset. Then $\CB^{\Omega}$ is indecomposable as a $\CZ$-module.
\end{lemma}

\begin{proof} We prove the theorem by induction on the number of elements in $\Omega$. If $\Omega=\{w\}$, then $\CB^{\Omega}=\CB^{w}=S$ is indecomposable. So let $\Omega\subset\CG$ be an arbitrary open set and let $x\in\Omega$ be minimal. Suppose the lemma is proven for $\Omega^\prime=\Omega\setminus\{x\}$. Every decomposition $\CB^{\Omega}=X\oplus Y$ induces a decomposition $\CB^{\Omega^\prime}=X^{\Omega^\prime}\oplus Y^{\Omega^{\prime}}$. By our inductive assumption we have, say, $X^{\Omega^\prime}=0$, hence $X=X^x$ is supported on $x$, and $\rho_{x,\delta x}|_{X^x}\colon X^x\to \CB^{\delta x}$ is the zero map.  So $\rho_{x,\delta x}|_{Y^x}\colon Y^x\to \CB^{\delta x}$ is surjective. Since $\rho_{x,\delta x}\colon \CB^x\to\CB^{\delta x}$ is  a projective cover, we have $X^x=0$, hence $X=0$ and  $\CB^{\Omega}$ is indecomposable. 
\end{proof}

\begin{lemma}\label{mult of B} For each $x\in\CV$ we have $\CB^x\cong S\oplus  \bigoplus_i S\{k_i\}$ for some $k_i<0$. If $\CB^x\cong S$, then $\CB^y\cong S$ for all $y\geq x$. 
\end{lemma}
\begin{proof} We claim that for each $x$ the $S$-modules $\CB^x$ and $\CB^{\delta x}$ are generated in non-negative degrees. This is certainly true for $\CB^w=S$. If the statement holds for $\CB^y$ for all $y>x$, then it follows first for $\CB^{\delta x}$ and then for its projective cover $\CB^x$. Hence for each $x$ we have $\CB^x\cong\bigoplus S\{k_i\}$ for some $k_i\leq 0$.

We have to show that there is exactly one copy of $S$ generated in degree zero, i.e.~ that $\CB^{x}_{\{0\}}=k$ for each $x$. Choose an enumeration $x_1,\dots,x_n$ of the vertices of $\CG$ such that $x_i<x_j$ implies $i<j$. Set $\Omega_i=\{x_j\mid j\geq i\}$.  We prove that $\CB^{x_i}_{\{0\}}=k$ and that $\Gamma(\Omega_i,\CB)_{\{0\}}=k$  by reverse induction on $i$. Since $\CB^{x_n}=\CB^w=S=\Gamma(\Omega_n,\CB)$, this is true for $i=n$. 

So let the claim be proven for $i$. Then for each open $\Omega\subset\Omega_i$ we have $\Gamma(\Omega,\CB)_{\{0\}}=k$, by the flabbiness of $\CB$. In particular, $\Gamma(\{> x_{i-1}\},\CB)_{\{0\}}=k$. Hence $\CB^{\delta x_{i-1}}_{\{0\}}=k$, so $\CB^{x_{i-1}}_{\{0\}}=k$, and $\Gamma(\Omega_{i-1},\CB)_{\{0\}}=k$ follows. Hence we proved the first part of the lemma.

Suppose that $\CB^x\cong S$. Then its quotient $\CB^{\delta x}$ is cyclic as an $S$-module. Choose $y\in\CV_{\delta x}$ and let $E\colon x\to y$ be the corresponding edge.  Then the canonical map $\CB^{\delta x}\to\CB^E$ (projection onto the $E$-coordinate) is surjective, hence $\CB^E$ is cyclic, so $\CB^y$ is cyclic since $\CB^E=\CB^y/ \alpha\cdot\CB^y$, where $\alpha$ is the label of $E$. Hence $\CB^y\cong S$. By induction we deduce that $\CB^y\cong S$ for all $y>x$.
\end{proof}

\begin{lemma}\label{smB} Suppose that $\CB^x\cong S$. Then $\CB_x=\alpha_1\cdots\alpha_{n}\cdot \CB^x$, where $\alpha_1,\dots,\alpha_{n}$ are the labels of the edges containing $x$. In particular, $\CB_x\cong \CM(x)\{-n\}$.
\end{lemma} 

\begin{proof}
If $\CB^x\cong S$, then $\CB^y\cong S$ for any $y>x$ by Lemma \ref{mult of B}, hence $\CB^E=S/\alpha_E\cdot S$ for any edge $E$ with vertex $x$, and $\CB^x\to\CB^E$ is isomorphic to the canonical map $S\to S/ \alpha(E)\cdot S$. Now $\CB_x\subset \CB^x$ is the intersection of the kernels of the maps $\rho_{x,E}\colon \CB^x\to\CB^E$ for all those edges $E$. Since their labels $\alpha_E$ are supposed to be pairwise linearly independent, $\CB_x=\alpha_1\cdots\alpha_{n}\cdot \CB^{x}$.
\end{proof}

\section{The ``smooth locus'' of $\CG$}\label{smloc}

Let $\CG$ be a finite connected GKM-graph with highest vertex $w$, and let $\CB$ be its intersection cohomology. Suppose that there is $l\in\DZ$ such that $\CB$ is self-dual of degree $l$, i.e.~ $\Bd\cong\CB\{l\}$. For a vertex $y\in\CV$ let  $n(y)\in\DN$  be the number of edges of $\CG$ containing $y$. The main result in this article is the following:

\begin{theorem}\label{mainT} Assume that $\CB$ is self-dual of degree $l$. For $x\in\CV$ the following are equivalent:
\begin{enumerate}
\item $\CB^x\cong S$.
\item $n(y)=l$ for all $y\geq x$. 
\end{enumerate}
\end{theorem}

\begin{proof} Let $\Omega=\{x\in\CV\mid n(y)=l\text{ for all $y\geq x$}\}$. Suppose that $\CB^x\cong S$. By  Lemma \ref{smB} we have $\CB_x\cong\CM(x)\{-n(x)\}$. The duals of the canonical maps $\CB\to\CB^x\cong\CM(x)$ and $\CM(x)\{-n(x)\}\cong\CB_x\to\CB$ are non-zero maps $\CM(x)\to\Bd$ and $\Bd\to\CM(x)\{n(x)\}$. They induce, by the self-duality of $\CB$, non-zero maps 
$$
\CM(x)\to\CB_x\{l\}\cong\CM(x)\{l-n(x)\}
$$
and 
$$
\CB^x\{l\}\cong\CM(x)\{l\}\to\CM(x)\{n(x)\}
$$
from which we deduce $0\leq l-n(x)$ and $l\leq n(x)$, respectively, hence $l=n(x)$. Since $\CB^x\cong S$ implies $\CB^y\cong S$ for all $y\geq x$, we deduce $x\in\Omega$. 

In order to prove the converse statement, let $1_\CB\in\CB$ be a generator of the homogeneous component of degree zero, and define $f\colon\CZ\to\CB$  by $f(z)=z.1_\CB$. Choose an isomorphism $\Bd\cong\CB\{l\}$ and consider the composition
$$
\phi\colon \CZ\{l\}\stackrel{f\{l\}}\longrightarrow\CB\{l\}\cong\Bd\stackrel{\dual f}\longrightarrow \Zd.
$$
\begin{lemma} Let $x\in\Omega$. Then $\phi^x\colon \CZ\{l\}^x\to(\Zd)^x$ is an isomorphism. Let $E\colon x\linie y$ be an edge with $x,y\in\Omega$. Then $\phi(E)\colon \CZ\{l\}(E)\to(\Zd)(E)$ is an isomorphism.
\end{lemma}
Before we prove the lemma, let us finish the proof of the theorem. The lemma implies that $\phi$ induces an isomorphism $\widehat{\CZ\{l\}^\Omega}\stackrel{\sim}\to\widehat{\Zd^{\Omega}}$ between the completions. But this map factors over $\CB\{l\}^\Omega=\widehat{\CB\{l\}^\Omega}$. Since the latter space is indecomposable by Lemma \ref{indB} we deduce $\CB^\Omega\cong\widehat{\CZ^\Omega}$ and, in particular, $\CB^x\cong S$ for all $x\in\Omega$.  
\end{proof}
\begin{proof}[Proof of the lemma]
The lemma follows easily from the following statements (2) and (3) by induction on the partial order on $\Omega$:
\begin{enumerate}
\item Suppose that $n(x)=l$ and $\phi^x\neq 0$. Then $\phi^x$ is an isomorphism.
\item For the highest vertex $w$, $\phi^w$ is an isomorphism.
\item Suppose that $\phi\colon x\linie y$ is an edge with $x,y\in \Omega$ and that $\phi^x$ is an isomorphism. Then $\phi^y$ and $\phi(E)$ are isomorphisms.  
\end{enumerate}

Ad (1): $\phi^x(1_\CZ)$ is a non-zero element in $(\Zd)^x$ of degree $-n(x)$. From Lemma \ref{strucDZ} we deduce $(\Zd)^x\cong\CM(x)\{l\}$. Hence $\phi^x$ is a non-zero endomorphism of $\CM(x)\{l\}$, hence must be an isomorphism.

Ad (2): The image of $1_\CZ$ in $\Bd$ under the composition $\CZ\{l\}\to\CB\{l\}\cong\Bd$ is an element of full support. Hence it does not vanish when restricted to $\CB_w=f(\CZ_w)$, so $\phi^w$ is non-zero. We have already deduced $n(w)=l$ from $\CB^w\cong S$. Hence $\phi^w$ is an isomorphism by (1). 

Ad (3): Let $u\in\CV$ be an arbitrary vertex and define $\xi^u:=\alpha_1\cdots\alpha_{n(u)}$, where $\alpha_1,\dots,\alpha_{n(u)}$ are the labels of all edges containing $u$. For any edge $E\colon x\stackrel{\alpha}\llinie y$ of $\CG$ define the element $z^E\in\CZ$ by $z^E_v=0$ for $v\not\in\{x,y\}$ and $z^E_x=z^E_y=\xi^x\cdot\xi^y\cdot\alpha^{-2}$ (this is the product of the labels of all edges containing $x$ or $y$ other than $E$). We will deduce statement (3) by evaluating elements of $\Zd$ on $z^E$. It might be helpful for the reader to take a look at the identity $(\ast)$ in Section \ref{dualofZ}.

Assume that $E\colon x\stackrel{\alpha}\llinie y$ is an edge with $x,y\in\Omega$ and that $\phi^x$ is an isomorphism. From $n(x)=l$ and Lemma \ref{strucDZ} we deduce that $\phi^x(1_\CZ)$ is a non-zero scalar multiple of $(\xi^x)^{-1}$. Suppose that $\phi^y=0$ . Then $\phi(1_\CZ)(z^E)$ is a non-zero multiple of $\xi^y\cdot \alpha^{-2}$. But it must be an element in $S$, which leads to a contradiction as $\xi^y$ is not divisible by $\alpha^2$. Hence $\phi^y\neq 0$ and $\phi^y$ is an isomorphism by (1).

Using the isomorphisms $\phi^x$ and $\phi^y$ we can identify $\phi(E)$ with an inclusion
$$
\CZ(E)\{l\}\subset  (\Zd)(E)\subset S\{l\}\oplus S\{l\}.
$$
In order to show that $\CZ(E)\{l\}=(\Zd)(E)$ it is sufficient to show that if $(\lambda,0)\in (\Zd)^{x,y}$, then $\lambda$ is divisible by $\alpha$ in $S$. Evaluating a preimage of $(\lambda,0)$ in $\Zd$ on $z^E$ shows that $\lambda\cdot \xi^y\cdot \alpha^{-2}\in S$. Since $\xi^y$ is divisible by $\alpha$, but not by $\alpha^2$, $\lambda$ must be divisible by $\alpha$, as was to be shown. 
\end{proof}

The next section contains the main  application of Theorem \ref{mainT}.

\section{Lusztig's conjecture}\label{Lconj}

In this section we apply Theorem \ref{mainT} to Lusztig's conjecture. References for the following notions and constructions are \cite{ModRep} and \cite{CharIrr} .

\subsection{Affine root systems}\label{reffaith}
Let $V$ be a finite dimensional rational vector space, and let $V^\ast:=\Hom_\DQ(V,\DQ)$ be the dual space. We denote by $\langle\cdot,\cdot\rangle\colon V\times V^\ast\to\DQ$ the natural pairing. 
Let $R\subset V$ be an irreducible and reduced root system. For $\alpha\in R$ we denote by $\alpha^\vee\in V^\ast$ the corresponding coroot. The  weight and coweight lattices are 
\begin{align*}
X&:=\{\lambda\in V\mid \langle\lambda,\alpha^\vee\rangle\in\DZ \text{ for all $\alpha\in R$}\},\\
X^\vee&:=\{v\in V^\ast\mid \langle\alpha, v\rangle\in\DZ \text{ for all $\alpha\in R$}\}.
\end{align*}

For $\alpha\in R$ and $n\in \DZ$ define the affine transformation $s_{\alpha,n}\colon V\to V$ by
$s_{\alpha,n}(\lambda)=\lambda-(\langle\lambda,\alpha^\vee\rangle -n)\alpha$. The {\em affine Weyl group} is the subgroup $\hCW$ of affine transformations on $V$  generated by the $s_{\alpha,n}$ with $\alpha\in R$ and $n\in \DZ$. 

We are now going to linearize the above affine action of $\hCW$ on
$V$. So we set $\hV:=V\oplus \DQ$ and define a linear action of $s_{\alpha,n}$ on $\hV$ by
$$
s_{\alpha,n}(\lambda,\xi)=(\lambda-(\langle\lambda,\alpha^\vee\rangle-\xi n)\alpha,\xi).
$$
This extends to a linear action of $\hCW$ on $\hV$ that leaves the level spaces $\hV_\kappa:=\{(\lambda,\kappa)\mid \lambda\in V\}$ for $\kappa\in\DQ$ stable. On $\hV_1$ we recover the affine action of $\hCW$ we started with. 

The hyperplane in $\hV$ fixed by $s_{\alpha,n}$ is 
$$
\hH_{\alpha,n}=\{(\lambda,\xi)\in \hV\mid \langle\lambda,\alpha^\vee\rangle=\xi n\}.
$$
We identify the dual space $\hV^\ast$ of $\hV$ with  $V^\ast\oplus\DQ$
in the obvious way and we consider $V^\ast$ as a subspace in  $\hV^\ast$ by the embedding of the first summand. If we define $\delta^\vee:=(0,1)\in \hV^\ast$, then 
$$
\alpha^\vee_n:=\alpha^\vee-n\delta^\vee\in\hV^\ast
$$
is a linear equation for $\hH_{\alpha,n}$. We call $\alpha^\vee_n$ the affine coroot associated to $\alpha$ and $n$. It is contained in the affine coweight lattice $\hX^\vee:=X^\vee\oplus \DZ\subset \hV^\ast$.

\subsection{The associated moment graph} Let $k$ be an arbitrary field
and let $W_k:=\hX^\vee\otimes_\DZ k$ be the $k$-vector space
associated to the lattice $\hX^\vee$.  For $v\in \hX^\vee$ we denote
by $\ol v$ the element $v\otimes 1\in W_k$. To the affine root system $\hR$ and the field $k$ we now associate the following (infinite) moment graph $\hCG_k$ over $W_k$:
\begin{itemize}
\item Its set of vertices is the affine Weyl group $\hCW$.
\item The vertices $x$ and $y$ are connected by an edge if $x=s_{\alpha,n}y$ for some $\alpha_n\in \hR^+$. This edge is labelled by $\ol{\alpha_n^\vee}\in W_k$.
\item The order on the set of edges is the affine Bruhat order on $\hCW$.
\end{itemize}
(In order to define the affine Bruhat order we have to choose a system
$R^+\subset R$ of positive roots.)
For $w\in\hCW$ we define $\hCG_{k,\le w}$ as the full submoment graph
of $\hCG_k$ whose set of edges is $\{\le w\}$. It is finite and
contains a highest vertex. We will later assume that $(k,w)$ are such
that $\hCG_{k,\le w}$ is a GKM-graph so that we are able to apply the
theory developed before. 

In particular, we then  have a Braden-MacPherson sheaf $\SB_{k,w}$ on
$\hCG_{k,w}$ and its space of global sections $\CB_{k,w}$. By construction, each stalk $\CB_{k,w}^x=\SB_{k,w}^x$ is a graded free $S=S(W_k)$-module of finite rank, i.e.~ there are $l_1,\dots, l_n\in\DZ$ such that $\CB_{k,w}^x\cong\bigoplus_{i=1,\dots, n} S\{l_i\}$. We define the {\em graded rank} of $\CB_{k,w}^x$ as
$$
\grk\,\CB_{k,w}^x:=\sum_{i=1,\dots,n} v^{2l_i}.
$$
For $x,w\in\hCW$ let $P_{x,w}\in\DZ[v]$ be the Kazhdan--Lusztig polynomial defined in \cite{KL1}. Here is our main conjecture:

\begin{conjecture}\label{conj-MC} Suppose that $k$ and $w$ are such that $\hCG_{k,w}$ is a GKM-graph. Then we have
$$
\grk\, \CB_{k,w}^x=P_{x,w}(v^{-2})
$$
for all $x\le w$.
\end{conjecture}

By construction, the degree of $P_{x,w}$  is strictly smaller than $1/2 (l(w)-l(x))$ if $x\ne w$, and $P_{w,w}=1$.  The polynomial $h_{x,w}\in \DZ[v]$ that is sometimes used instead of $P_{x,w}$ in the formulation of the above conjecture satsifies  
$h_{x,w}(v)=v^{l(w)-l(x)}P_{x,w}(v^{-2})$.

For $x\le w$ we define $n_w(x)$ as the number of edges of $\hCG_{k,\le w}$ that contain $x$ as a vertex. By definition, we have
$$
n_w(x)=\{t\in\hCT\mid tx\le w\},
$$
where we denote by $\hCT\subset \hCW$ the reflections, i.e.~ the set
$\{s_{\alpha,n}\mid \alpha\in R,n\in\DZ\}$. Let us denote by $l\colon
\hCW\to\DN$ the affine length function. The following is a conjecture
of Deodhar, which is proven in \cite{C}.

\begin{theorem}\label{theorem-conjDeo} For each $x\le w$ we have
$$
P_{x,w}=1\text{ if and only if } n_w(x)=l(w).
$$
\end{theorem}
The above equivalence cannot yet be shown for arbitrary Coxeter systems (the proof in \cite{C} uses the positivity of the coefficients of $P_{x,w}$).
Note that $P_{x,w}=1$ if and only if $P_{x,w}(1)=1$. In order to apply
our main Theorem  \ref{mainT} we need yet another  result:

\begin{theorem}  \cite[Theorem 6.1]{Fie2} Suppose that $k$ and $w$ are such that $\hCG_{k,w}$ is a GKM-graph. Then 
$\CB_{k,w}$ is self-dual of degree $l(w)$. 
\end{theorem}
(Note that the module $B(w)$ in \cite{Fie2} is $\CB_{k,w}\{l(w)/2\}$ by
Corollary 6.5 in {\em loc.cit} and our grading normation $\deg W_k=1$.)
Hence we can apply Theorem \ref{mainT} in our situation. Together with Theorem \ref{theorem-conjDeo} we deduce the multiplicity one case of Conjecture \ref{conj-MC}:

\begin{theorem}\label{theorem-multone} Suppose that $k$ and $w$ are such that $\hCG_{k,w}$ is a GKM-graph. For all $x\le w$ we have
$$
\grk\,\CB_{k,w}^x(1)=1 \text{ if and only if } P_{x,w}(1)=1.
$$
\end{theorem}

Finally, we apply the above result to Lusztig's conjecture.

\subsection{Jordan--H\"older multiplicities of baby Verma modules}
Let $\fg$ be the Lie algebra over $k$ with root system $R$. For
$\lambda\in X$ let $\Delta(\lambda)$ be the ($X$-graded) baby Verma
module with highest weight $\lambda$ and denote by $L(\lambda)$ its
simple quotient. Let $\rho=1/2\sum_{\alpha\in R^+}\alpha\in X$ be the half-sum of positive roots. By $[\Delta(\lambda):L(\mu)]$ we denote the Jordan-H\"older multiplicity of $L(\mu)$ in $\Delta(\lambda)$. The multiplicity version of Lusztig's conjecture on the irreducible rational characters of the simply connected, connected algebraic group over $k$ with root system $R$  is the following:
\begin{conjecture} Suppose that the characteristic of $k$ is bigger than the Coxeter number of $R$. For $x,w\in\hCW$ we have
$$
[\Delta(x(\rho)-\rho):L(w(\rho)-\rho)]=p_{x,w}(1).
$$
\end{conjecture}
Here $p_{x,w}$ is the periodic polynomial defined by Lusztig (\cite{L2}, cf.~  \cite{CharIrr} for notation and normalization). Due to an inherent symmetry in the theory it is enough to prove the above conjecture for $w\in\hCW^{\res,-}$, where 
$$
\hCW^{\res,-}:=\left\{w_0w\in\hCW\mid 0\le\langle w\cdot 0 , \alpha^\vee\rangle<
\ch\, k\text{ for all simple roots $\alpha$}\right\}
$$
is the set of restricted antidominant elements ($w_0$ is the longest
element in the finite Weyl group). The main result of \cite{ModRep} is the following:

\begin{theorem}\label{theorem-MTModRep} Suppose that the
  characteristic of $k$ is bigger than the Coxeter number. For $w\in\hCW^{\res,-}$ and $x\in \hCW$ we have
$$
[\Delta(x\cdot 0):L(w\cdot 0)]=\grk\, \CB_{k,w}^x(1).
$$
\end{theorem}

For $w\in\hCW^{\res,-}$ and $x\in\hCW$ we have 
$$
P_{x,w}(1)=p_{x,w}(1),
$$
even though $P_{x,w}\ne p_{x,w}$ in general (cf.~\cite{L2}).  This
curiosity is mirrored in the fact that the functor $\Phi$,
constructed in
\cite{ModRep}, preserves multiplicities but not the
gradings. Moreover, if the characteristic of $k$ is bigger than the
Coxeter number, then $\hCG_{k,w}$ is a GKM-graph for all
$w\in\hCW^{\res,-}$, cf.~\cite[Lemma 9.1]{ModRep}. So
we deduce from Theorem \ref{theorem-multone} and  Theorem \ref{theorem-MTModRep}  the multiplicity one case of Lusztig's conjecture:

\begin{theorem} Suppose that the characteristic of $k$ is bigger than the Coxeter number.  Then
$$
[\Delta(x\cdot 0):L(w\cdot 0)]=1\text{ if and only if } p_{x,w}(1)=1.
$$
\end{theorem}

\end{document}